\documentclass[a4paper,12pt]{amsart}
\usepackage{amsmath,amssymb,latexsym,amsfonts,amscd}

\title{Automorphisms of K3 surfaces}
\author{De-Qi Zhang}
\address{Department of Mathematics, National University of Singapore,
2 Science Drive 2, Singapore 117543, Singapore}
\email{matzdq@math.nus.edu.sg}

\thanks{
The author was partially supported by an Academic
Research Fund of NUS}
\newcommand\Pic{\text{\rm Pic}}

\newcommand\OO{{\mathcal{O}}}

\newcommand\ZZ{{\mathbb{Z}}}

\newcommand\rk{{\text{\rm rank}}}
\newcommand\diag{{\text{\rm diag}}}
\newcommand\Sing{{\text{\rm Sing}}}
\newcommand\Tr{{\text{\rm Tr}}}
\newcommand\Aut{{\text{\rm Aut}}}
\newcommand\Gal{{\text{\rm Gal}}}
\newcommand\Bir{{\text{\rm Bir}}}
\newcommand\Ker{{\text{\rm Ker}}}
\newcommand\ord{{\text{\rm ord}}}
\newcommand\Imm{{\text{\rm Im}}}

\newtheorem{thm}{Theorem}[section]
\newtheorem{lem}[thm]{Lemma}

\newtheorem{prop}[thm]{Proposition}

\theoremstyle{definition}

\newtheorem{prob}[thm]{Problem}
\newtheorem{rem}[thm]{Remark}
\theoremstyle{remark}

\newtheorem{ex}[thm]{Example}
\newtheorem{definition}[thm]{Definition}
\begin{document}
\begin{abstract}
In this note, we report some progress we made recently on
the automorphisms groups of $K3$ surfaces. A short and straightforward 
proof of the impossibility of
$\ZZ/(60)$ acting purely non-symplectically on a $K3$ surface, is also given,
by using Lefschetz fixed point formula for vector bundles.

\end{abstract}
%%%%%%%%%%%%%%%%%%%%
\maketitle
%%%%%%%%%
\pagestyle{myheadings} \markboth{\hfill D. -Q.
Zhang\hfill}{\hfill Automorphisms groups of $K3$ surfaces\hfill}
%%%%%%%%%%%%%%%%%%%%%%%%%%%%%%%%%%

\section{\bf Introduction}

We work over the complex numbers field ${\bold C}$.
A smooth projective surface $X$ is a $K3$ surface if
the canonical line bundle (is trivial) $\omega_X \cong \OO_X$ and
if the irregularity $h^1(X, \OO_X) = 0$. In this note
we shall report some progress we made recently (mostly with K. Oguiso)
on automorphisms of $K3$ surfaces.

\par
By abuse of notation, we denote by the same $\omega_X$ a nowhere vanishing
global holomorphic $2$-form on $X$. So $H^0(X, {\omega_X}) = {\bold C} \omega_X$.
Now the Hodge decomposition says that $H^2(X, {\bold C}) = {\bold C} \omega_X \oplus
H^{1,1}(X) \oplus {\bold C} \overline{\omega_X}$.

\par
Let $\Aut(X)$ be the group of all automorphisms of a $K3$ surface $X$.
Take a subgroup $G \subseteq \Aut(X)$. For every $g$ in $G$, we have
$g^* \omega_X = \alpha(g) \omega_X$ for some $\alpha(g)$ in the multiplicative group
${\bold C}^* = {\bold C} - \{0\}$. One sees that $\alpha : G \rightarrow {\bold C}^*$ is 
a group homomorphism. Let $G_N = \Ker(\alpha)$ be the kernel.
By Nikulin \cite{Ni-finite}, Theorem 0.1 or Sterk \cite{St}, Lemma 2.1, 
the image $\Imm (\alpha)$ is a finite
subgroup (of order $I = I(G) = I(X, G)$) in ${\bold C}^*$, so it equals the multiplicative group
$\mu_I = \langle \zeta_I \rangle$,
where $\zeta_I = \exp( 2 \pi \sqrt{-1}/I)$ is an $I$-th primitive root of $1$.
Then one has the {\bf basic exact sequence} of groups:
$$1 \rightarrow G_N \rightarrow G \rightarrow \mu_I \rightarrow 1.$$

One may say that the study of $G$ is reduced to that of $G_N$ and $\mu_I$,
where $G_N$ acts on $X$ symplectically (see Notation and Terminology below)
while $h^* \omega_X = \zeta_I \omega_X$
for some pre-image $h$ in $G$ of a generator in $\mu_I$.

\begin{prob}\label{}[{\bf Extension Problem}]
Given symplectic subgroup $H$ and a cyclic subgroup $\mu_I$ of $\Aut(X)$ on a $K3$ surface $X$,
find the conditions on $H$ and $I$ such that there is a subgroup $G$ of $\Aut(X)$
with $G_N = H$ and $I = I(X, G)$.
\end{prob}

In the forthcoming paper Ivanov-Oguiso-Zhang \cite{IOZ}, we will tackle this problem.

Now a summary of the contents. We start with symplectic actions on $K3$ surfaces,
then purely non-symplectic actions, and actions of mixed type. Finally
we mention the connection between $K3$ surfaces with a non-symplectic action
and log Enriques surfaces initially defined in Zhang \cite{Z-enriq} and \cite{Z-enriq2}.
Examples of actions are either included or mentioned about their whereabouts.
In \S 3, we also provide a new proof of the impossibility of
$\ZZ/(60)$ acting purely non-symplectically on a $K3$ surface,
by using Lefschetz fixed point formula for vector bundles.

\vskip 1pc \noindent
{\bf Notation and Terminology}

\par
For a $K3$ surface $X$, we use the same symbol
$\omega_X$ ($\cong \OO_X$) to denote the canonical line bundle and also
a generator of $H^0(X, \omega_X) = {\bold C} \omega_X$.
Let $S_X = \Pic(X) < H^2(X, {\bold Z})$ be the Neron Severi lattice
and $T_X = (S_X)^{\perp} = \{t \in H^2(X, {\bold Z}) \, | \, (t, s) = 0, \, \forall s \in S_X \}$
the transcendental lattice.

\par
A group $G$ is called a $K3$ group if $G$ is 
a subgroup of $\Aut(X)$ for some $K3$ surface $X$.
An element $g$ of $\Aut(X)$ is simplectic if
$g^* \omega_X = \omega_X$. A group $G$ is symplectic if
every element of $G$ is symplectic.
An element $g$ (or the group $\langle g \rangle$) of order $n$ is purely non-symplectic
if $g^* \omega_X = \eta_n \omega_X$ for some primitive $n$-th root of $1$.

\par
For an integer $I \ge 2$, the Euler function
$\varphi(I) = \#\{n \, | \, 1 \le n \le I - 1; \gcd(n, I) = 1\}$.
Set $\zeta_I = \exp(2 \pi \sqrt{-1}/I)$, a primitive $I$-th root of $1$.
Set $\mu_I = \langle \zeta_I \rangle$, a cyclic group of order $I$.

\par
For groups $H$ and $K$, the symbol $H . K$ denotes a group $G$
such that $H$ is normal in $G$ and $G/H \cong K$, i.e., $G$ is an extension
of $K$ by $H$.
When $K$ can be chosen to be a subgroup of $G$, then we denote 
$H : K$, i.e., $G$ is a split extension of $K$ by $H$.

\par
$2^n$, $3^n$, etc denote elementary abelian groups $({\bold Z}/(2))^{\oplus n}$, 
$({\bold Z}/(3))^{\oplus n}$, etc.

\par
$S_n$ is the symmetric group (of order $n!$) in $n$ letters.
$A_n$ is the alternating group (of order $n!/2$) in $n$ letters.
$M_{24}, M_{23}, M_{20}, M_{10}$, etc are Mathieu groups.
$Q_8$ is the quaternion group of order 8, while
$D_n$ is the dihedral group of order $n$ (see Mukai \cite{Mu-sym}).

\par \vskip 1pc \noindent
{\bf Acknowledgement.} 
\par \noindent
The author would like to thank Professors S. Kondo and K. Oguiso
for the valuable comments.

\section{\bf Symplectic $K3$ groups}

The classification of abelian symplectic groups was done by Nikulin \cite{Ni-finite}.
Especially, he proved:

\begin{thm}\label{Ni-fix}
Let $g$ be a symplectic automorphism on a $K3$ surface $X$.
Then $\ord(g) \le 8$. Moreover, the fixed locus (point wise) $X^g$ 
has cardinality $8, 6, 4, 4, 2, 3$, or $2$, if $\ord(g)$ equals
$2, 3, 4, 5, 6, 7$, or $8$, respectively.
\end{thm}

\par
Mukai \cite{Mu-sym} has determined all maximum finite symplectic $K3$ groups (11 of them);
see also Kondo \cite{Ko-lat} for a lattice-theoretic proof:

\begin{thm}\label{Mu-sym}
For a finite group $G$, the following are equivalent:

\par \noindent
$(1)$ $G$ has a (faithful) symplectic action on a $K3$ surface.

\par \noindent
$(2)$ $G$ has an embedding $G \subset M_{23}$ into the Mathieu group
and decomposes $\{1, 2, 3, \dots, 24\}$ 
into at least $5$ orbits.
(The bigger Mathieu group $M_{24}$
acts naturally on the $24$-element set with the smaller $M_{23}$ as a one-point stabilizer subgroup.)

\par \noindent
$(3)$ $G$ has a Mathieu representation $V$ over {\bf Q} with $\dim V^G \ge 5$ and
the $2$-Sylow subgroups of $G$ can be embedded into $M_{23}$.

\par \noindent
$(4)$ $G$ is isomorphic to one of the $11$ (maximum symplectic) groups below
where the order is indicated in the parenthesis:

\par \vskip 1pc
$\text{\rm PSL}_2(7)$ $(168)$, \,\, $A_6$ $(360)$, \,\, $S_5$ $(120)$, \,\, 
$M_{20} = 2^4 : A_5$ $(960)$,

\par \vskip 1pc
$F_{384} = 2^4 : S_4$ $(384)$, \,\,
$A_{4, 4} = (S_4 \times S_4) \cap A_{8}$ $(288)$,

\par \vskip 1pc
$T_{192} = (Q_8 * Q_8) : S_3$ $(192)$, \,\, 
$H_{192} = 2^4 : D_{12}$ $(192)$,

\par \vskip 1pc
$N_{72} = 3^2 : D_8$ $(72)$, \,\, 
$M_9 = 3^2 : Q_8$ $(72)$, \,\,
$T_{48} = Q_8 : S_3$ $(48)$.
\end{thm}

\par
The result above completely classifies all finite symplectic actions
on $K3$ surfaces. The possible subgroups $G$ (exactly 80 of them) of the 11 maximum symplectic
$K3$ groups acting on some $K3$ surface $X$ together with the data
of $\Sing(X/G)$ has been given by Xiao \cite{Xi1}, the list.

\section{\bf Purely non-symplectic $K3$ groups}

We now turn to the purely non-symplectic $K3$ groups.
Suppose that a cyclic group $\mu_I = \langle g \rangle$ ($I \ge 2$) acts on a $K3$ surface
$X$ such that $g^* \omega_X = \eta_I \omega_X$, where 
$\eta_I$ is a primitive $I$-th root of 1.
By Nikulin \cite{Ni-finite} Theorem 0.1, or Sterk \cite{St} Lemma 2.1, 
the Euler function $\varphi(I) \, | \, \rk(T_X)$.
Since $\rk(T_X) = 22 - \rk(S_X) \le 21$, 
we get:

\begin{lem}\label{index}
Suppose $g$ is an order-$I$ ($I \ge 2$) automorphism of a $K3$ surface $X$
such that $g^* \omega_X = \eta_I \omega_X$ with $\eta_I$ a primitive $I$-th
root of $1$. Then $\varphi(I) \le 21$ and hence $I$ equals one of the following:

\par \noindent
Case $\varphi(I) = 20$. $I = 66, 50, 44, 33, 25$.

\par \noindent
Case $\varphi(I) = 18$. $I = 54, 38, 27, 19$.

\par \noindent
Case $\varphi(I) = 16$. $I = 60, 48, 40, 34,32, 17$.

\par \noindent
Case $\varphi(I) = 12$. $I = 42, 36, 28, 26, 21, 13$.

\par \noindent
Case $\varphi(I) = 10$. $I = 22, 11$.

\par \noindent
Case $\varphi(I) = 8$. $I = 30, 24, 20, 16, 15$.

\par \noindent
Case $\varphi(I) = 6$. $I = 18, 14, 9, 7$.

\par \noindent
Case $\varphi(I) = 4$. $I = 12, 10, 8, 5$.

\par \noindent
Case $\varphi(I) = 2$. $I = 6, 4, 3$.

\par \noindent
Case $\varphi(I) = 1$. $I = 2$.
\end{lem}

Conversely, each $I$ with $I \ne 60$ in the lemma above is geometrically realizable
(see Kondo \cite{Ko92}, or Machida-Oguiso \cite{MO}, Proposition 4]).

\begin{prop}
For each $I$ with $I \ne 60$ in the lemma above, there is a $K3$ surface $X$
such that $\Aut(X) \ge \mu_I = \langle g \rangle$
and $g^* \omega_X = \zeta_I \omega_X$, where 
$\zeta_I = \exp(2 \pi \sqrt{-1}/I)$ is a primitive $I$-th root of 1.
\end{prop}

Indeed, $\mu_{60}$ can not act purely non-symplectically by
the result below (the early proofs are Machida-Oguiso \cite{MO}, Theorem 5.1]
and Xiao \cite{Xi2}, Theorem, page 1); our proof here comes from the
very simple idea: the calculation of Lefschetz number by Atiyah, Segal and Singer,
\cite{AS1}, \cite{AS2},
for the trivial line bundle, and hence more straightforward.

\par
Therefore, the lemma and proposition above give the
complete classification for the $\mu_I$ appearing in the basic exact sequence in
the Introduction.

\begin{prop}\label{I60} (see \cite{MO} and \cite{Xi2}).
There is no $K3$ surface with a purely non-symplectic action
by the group $\mu_{60}$. Namely, there is no pair $(X, \langle g \rangle)$
of a $K3$ surface $X$ and an order $60$ automorphism $g$ with
$g^* \omega_X = \zeta_{60} \omega_X$.
\end{prop}

\begin{proof}
Suppose the contrary that there is a pair $(X, \langle g \rangle)$ as in the statement.
Set $\omega = \omega_X$ and $\zeta = \zeta_{60}$, so $g^* \omega = \zeta \omega$.
Let $\{C_i\}$ be the set of all $g$-fixed (point wise) irreducible curves.
So at a point $Q$ on $C_i$, one can diagonalize $g^*$ as $(1, \zeta)$
with appropriate coordinates.
Let $n = \sum (1 - g(C_i))$;
we set $n = 0$ when there is no such $g$-fixed curve $C_i$.
Let $M_j = \{P_{jk}\}$ be the set of points $P$ at which $g^*$ can be diagonalized as
$(\zeta^{-j}, \zeta^{j+1})$ with $1 \le j < 60/2 = 30$. Set $m_j = |M_j|$ and $m = \sum m_j$.
Then the fixed locus $X^g$ (point wise) is the disjoint union of smooth curves $C_i$
and $m$ isolated points $P_{jk}$ 
(see, e.g., Oguiso - Zhang \cite{OZ-AJM}, Lemma 2.1 (3) for the smoothness).
By Atiyah - Segal \cite{AS1} and Atiyah - Singer \cite{AS2} at pages 542 and 567, 
we can calculate the Lefschetz number
$L(g)$ in two different ways:
$$L(g) = \sum_{i=0}^2 (-1)^i \, \Tr(g^* | H^i(X, \OO_X)),$$
$$L(g) = \sum_{j, k} a(P_{jk}) + \sum_i b(C_i),$$
$$a(P) = \frac{1}{\det(1 - g^* | T_{P})},$$
$$b(C_i) = \frac{1 - g(C_i)}{1 - \zeta} - \frac{\zeta C_i^2}{(1 - \zeta)^2} 
= \frac{(1 - g(C_i))(1 + \zeta)}{(1 - \zeta)^2}.$$

Here $T_P = T_{X, P}$ is the tangent space at $P$,
and $\zeta^{-1}$ is the eigenvalue of the action $g_*$ on the normal bundle of $C_i$.
We see then that $\sum_i b(C_i) = n (1 + \zeta)/(1 - \zeta)^2$.
Also $\sum a(P_{jk}) = \sum_j m_j/(1 - \zeta^{-j})(1 - \zeta^{j+1})$.
On the other hand, by the Serre duality $H^2(X, \OO_X) = H^0(X, \omega)^{\vee}$, we have 
$L(g) = 1 + \zeta^{-1}$.
Identify the two $L(g)$, so $x := \zeta$ satisfies the equation below:
$$F_1(x) = -(1 + x^{-1}) + 
\sum_{j=1}^{29} \frac{m_j}{(1 - x^{-j})(1 - x^{j+1})} + \frac{n(1+x)}{(1-x)^2} = 0.$$

Note that the minimal polynomial $\Phi_{60}(x)$ of $x = \zeta$ over {\bf Q} is given as follows;
this is also obtained by factorizing $x^{30} + 1$:
$$\Phi_{60}(x) = x^{16} + x^{14} - x^{10} - x^8 - x^6 + x^2 + 1.$$
Our $x = \zeta$ also satisfies the following equations (factorizing $(x^{10})^3 + 1$):
$$x^{30} + 1 = 0 = x^{20} - x^{10} + 1.$$
Since $x^{\pm 30} = -1$, we have $(1 - x^{-i})(1 - x^{i+1}) = (1 + x^{30-i})(1 + x^{i+1-30})$.
Substituting these equalities for $15 \le i \le 29$ into the equation $F_1(x) = 0$,
we obtain a new equation
$0 = F_1(x) = F_2(x)/G_2(x)$, where $F_2(x)$ is a polynomial of degree $146$ in $x$
with coefficients in $m_j$, and $G_2(x)$ is the product of polynomials each of
which is of degree $\le 12$ in $x$. Since the minimal polynomial $\Phi_{60}(x)$ of
$x = \zeta$ is of degree $16$, our $x = \zeta$ satisfies $G_2(x) \ne 0$ and 
$F_2(x) = 0$. Substituting $x^{30 k} = 1$ (resp. $-1$) when $k$ is even (resp. odd)
into the equation $F_2(x) = 0$, we obtain a polynomial equation
$F_3(x) = 0$ of degree $\le 29$ in $x$. Substituting $x^{20} = x^{10} - 1$
into the equation $F_3(x) = 0$, we obtain an equation $F_4(x) = 0$
of degree $\le 19$ in $x$. Finally substituting $\Phi_{60}(x) = 0$ into
the equation $F_4(x) = 0$, we obtain an equation
$F_5(x) = \sum_{i = 0}^{15} d_i x^i = 0$ of degree $\le 15$ in $x$,
where each $d_i$ is an integer linear combination of $m_j$'s.
Note that $x^i = \zeta^i$ ($0 \le i \le 15$) are linearly independent over
{\bf Q}. So we get linear equations in $m_j$:
$d_0 = d_1 = \dots = d_{15} = 0$. Solving them, we arrive at
the following relations:
$$
(*) \hskip 0.5pc 4m_{1} = 
      -1 + 2 m_{2} +2 m_{20} - 2 m_{21} +4 m_{22} - 4 m_{23} +2 m_{26}$$
$$ - 2 m_{27} +4 m_{28} - 4 m_{29} - 2 m_{3} +4 m_{6} - 4 m_{7} +2 m_{8} - 
            2 m_{9} +8 n,$$

\par \vskip 0.5pc

$$  9m_{10} = 
      60 - 12 m_{17} +18 m_{18} - 17 m_{19} - 141 m_{2} - 94 m_{20} +
            81 m_{21} - 189 m_{22} $$
$$+198 m_{23} +6 m_{24} - 12 m_{25} - 81 m_{26} +
            81 m_{27} - 156 m_{28} +156 m_{29} $$
$$+63 m_{3} - 6 m_{4} - 24 m_{5} - 
            198 m_{6} +189 m_{7} - 105 m_{8} +96 m_{9} +312 n,$$

\par \vskip 0.5pc

$$ 36m_{11} = 
      267 - 48 m_{17} +108 m_{18} - 104 m_{19} - 582 m_{2} - 394 m_{20} +
            342 m_{21} - 792 m_{22} $$
$$+828 m_{23} +24 m_{24} - 48 m_{25} - 342 m_{26} +
            342 m_{27} - 624 m_{28} +624 m_{29} $$
$$+270 m_{3} - 24 m_{4} - 96 m_{5} - 
            828 m_{6} +792 m_{7} - 438 m_{8} +402 m_{9} +1248 n,$$

\par \vskip 0.5pc

$$  36m_{12} = 
      57 + 12 m_{17} - 72 m_{18} +56 m_{19} - 282 m_{2} - 134 m_{20} +
            162 m_{21} - 252 m_{22}$$
$$ +216 m_{23} - 24 m_{24} +12 m_{25} - 162 m_{26} +
            162 m_{27} - 384 m_{28} +384 m_{29}$$
$$ +90 m_{3} - 84 m_{4} +24 m_{5} - 
            288 m_{6} +252 m_{7} - 138 m_{8} +102 m_{9} +768 n,$$

\par \vskip 0.5pc

$$ 72m_{13} = 
     39 + 48 m_{17} - 144 m_{18} +104 m_{19} - 426 m_{2} - 182 m_{20} +
            234 m_{21} - 252 m_{22}$$
$$ +180 m_{23} - 60 m_{24} +84 m_{25} - 270 m_{26} +
            234 m_{27} - 600 m_{28} +600 m_{29}$$
$$ +126 m_{3} - 156 m_{4} +60 m_{5} - 
            324 m_{6} +252 m_{7} - 138 m_{8} +102 m_{9} +1200 n,$$

\par \vskip 0.5pc

$$ 18m_{14} = 
      -15 + 12 m_{17} - 36 m_{18} +32 m_{19} - 48 m_{2} - 2 m_{20} +$$
$$            18 m_{21} - 18 m_{23} - 24 m_{24} +30 m_{25} - 36 m_{26} +36 m_{27} - 
            96 m_{28}$$
$$ +87 m_{29} - 30 m_{4} +24 m_{5} +18 m_{6} - 36 m_{7} +6 m_{8} - 
            6 m_{9} +174 n, $$

\par \vskip 0.5pc

$$ 18m_{15} = 
      30 + 12 m_{17} - 36 m_{18} +32 m_{19} - 66 m_{2} - 56 m_{20} +
            72 m_{21} - 144 m_{22}$$
$$ +126 m_{23} +12 m_{24} - 6 m_{25} - 54 m_{26} +
            54 m_{27} - 96 m_{28} +87 m_{29} $$
$$+18 m_{3} +6 m_{4} - 12 m_{5} - 126 m_{6} +
            108 m_{7} - 48 m_{8} +48 m_{9} +174 n, $$

\par \vskip 0.5pc

$$ 72m_{16} = 
      -39 + 96 m_{17} - 8 m_{19} +138 m_{2} +86 m_{20} - 90 m_{21} +
            180 m_{22}$$
$$ - 180 m_{23} - 12 m_{24} +60 m_{25} +54 m_{26} - 90 m_{27} +
            168 m_{28} - 168 m_{29}$$
$$ - 54 m_{3} +12 m_{4} +12 m_{5} +180 m_{6} - 
            180 m_{7} +138 m_{8} - 102 m_{9} - 336 n.$$

\par \vskip 1pc 
To double check that there is no any human error in the above substitutions,
we substitute the above relations back into the original equation and get
$F_1(x) = \Phi_{60}(x) F_6(x)/G_6(x)$, where $F_6(x)$ is a polynomial of
degree $130$ in $x$ and $G_6(x)$ is a product of polynomials each of which is of
degree $\le 12$ in $x$. So no error should have occurred during the substitution process,
noting that $x = \zeta$ satisfies its minimal polynomial (over {\bf Q}) $\Phi_{60}(x) = 0$.
The actual simple linear algebra calculation is done by the "Maple",
though a calculation by hand is possible if the reader is patient and careful enough
to avoid any error in doing the simple arithmetics and Gaussian elimination.

\par
Now to conclude the proof, we note that the first equality (*) in the final relations above says that
$4m_1 = - 1 + 2 s$, with $s$ an integer (a integer linear combination of $m_i$ and $n$). 
This is impossible because $m_1$ is an integer. We reach a contradiction.
So there is no such pair $(X, \mu_{60})$ and the proposition is proved.
\end{proof}

\begin{rem}
As pointed out by Oguiso (and done in \cite{MO}), 
indeed, one can say a bit stronger than the proposition above
(consider $X/G_N$):
If $X$ is a $K3$ surface and $G \le Aut(X)$ is finite then $\mu_{I(X,G)} \ne \mu_{60}$.
This is because one always has $G_N = (1)$ and $G = \mu_{I(X,G)}$ whenever $\varphi(I(X, G)) \ge 10$
and $I(X, G) \ne 28$; see Ivanov-Oguiso-Zhang \cite{IOZ}, or Machida-Oguiso \cite{MO}
the proof of Lemma (4.1).
\end{rem}

\par
For a $K3$ surface $X$, the cohomology $H^2(X, {\bold Z}) = U^{\oplus 3} \oplus E_8^{\oplus 2}$
is an even lattice, where $U = {\bold Z} e_1 + {\bold Z} e_2$ with $e_i^2 = 0$ and $e_1 . e_2 = 1$,
and $E_8$ is the negative definite even lattice of rank 8.
So the study of $\Aut(X)$ is reduced to that of $\Aut(S_X)$ to some extent.
Indeed, by Pjateckii-Shapiro and Shafarevich \cite{PS} (see also Burns - Rapoport \cite{BR}) 
up to finite
index and finite co-index, $\Aut(X)$ is isomrphic to $\Aut(S_X) / \Gamma(X)$
where $\Gamma(X)$ is generated by reflections associated to $(-2)$-vectors.

\par
Let us consider the natural exact sequence
$$(1) \rightarrow H_X \rightarrow \Aut(X) \rightarrow \Aut(S_X).$$ 

\par
Note that $H_X$ acts faithfully on $T_X$ and on the 1-dimensional space
${\bold C} \omega_X$ (see Nikulin \cite{Ni-finite} Cor. 3.3, or Sterk \cite{St} Lemma 2.1). 
Indeed, $H_X$
is a cyclic finite group of order $h_X$ with the Euler function
$\varphi(h_X) \, | \, \rk(T_X)$.
The results below show that these $X$ with maximum
$\varphi(h_X) = rk (T_X)$ are unique up to isomorphisms.
For the unique examples of the $X$'s in the results below,
we refer to Kondo \cite{Ko92}; see also Examples \ref{ex1} and \ref{ex2} below.

\begin{thm} (see Kondo \cite{Ko92}, Main Th).

\par
Set $\Sigma := \{66, 44, 42, 36, 28$, $12\}$. 

\par \noindent
$(1)$ Let $X$ be a $K3$ surface with $\varphi(h_X) = \rk (T_X)$
and with unimodular transcendental lattice $T_X$. Then $h_X \in \Sigma$.

\par \noindent
$(2)$ Conversely, for each $N \in \Sigma$, there exists, modulo isomorphisms, a
unique $K3$ surface $X$ such that $h_X = N$ and $\varphi(h_X)
= \rk(T_X)$. Moreover, $T_X$ is unimodular for this $X$.
\end{thm}

\begin{thm} (see Oguiso - Zhang \cite{OZ-PAMS}, Theorem $2$).

Set $\Omega := \{3^k (1 \le k \le 3), \, 5^{\ell} (\ell = 1,2), \, 7, 11, 13, 17, 19\}$. 

\par \noindent
$(1)$ Let $X$ be a $K3$ surface with $\varphi(h_X) = \rk (T_X)$
and with non-unimodular transcendental lattice $T_X$. Then $h_X \in \Omega$.

\par \noindent
$(2)$ Conversely, for each $N \in \Sigma$, there exists, modulo isomorphisms, a
unique $K3$ surface $X$ such that $h_X = N$ and $\varphi(h_X)
= \rk(T_X)$. Moreover, $T_X$ is non-unimodular for this $X$.
For $I = 19, 17$, or $13$, the unique pair is given in Example \ref{ex1} below.
\end{thm}

By Theorem \ref{Ni-fix} and Lemma \ref{index}, every prime order $p$ automorphism $g$
of a $K3$ surface $X$ satisfies $p \le 19$. Let us give examples for the
largest three $p$ (see Kondo \cite{Ko92}).

\begin{ex}\label{ex1}
Here are examples of pairs $(X_p, \mu_p)$ ($p = 19, 17$ or 13)
of the $K3$ surface $X_p$ with a $\mu_p = \langle g_p \rangle$ action, given by its Weierstrass equation and the action on the coordinates.

\par \vskip 1pc

$X_{19}: \,\, y^2 = x^3 + t^7 x + t$, \,\, 
$g_{19}^*(x, y, t) = (\zeta_{19}^7 x, \zeta_{19} y, \zeta_{19}^2 t)$.

$X_{17}: \,\, y^2 = x^3 + t^7 x + t^2$, \,\, 
$g_{17}^*(x, y, t) = (\zeta_{17}^7 x, \zeta_{17}^2y, \zeta_{17}^2 t)$.

$X_{13}: \,\, y^2 = x^3 + t^5 x + t^4$, \,\, 
$g_{13}^*(x, y, t) = (\zeta_{13}^5 x, \zeta_{13} y, \zeta_{13}^2 t)$.
\end{ex}

Indeed, these pairs are unique with no conditions on $h_X$ or so on;
but they turn out to satisfy the conditions in the above theorem:

\begin{thm} (Oguiso-Zhang \cite{OZ-PAMS} Cor. $3$).
Let $X$ be a $K3$ surface with an automorphism $g$ of order $p = 19$, $17$ or $13$
(the three largest possible prime orders. Then the pair $(X, \langle g \rangle)$
is isomorphic to the pair $(X_p, \langle g_p \rangle)$ in the above example.
\end{thm}

\begin{rem}
A similar result as above does not hold for $p = 11$. Indeed, there are three 1-dimensional
families of pairs $(X_t, \mu_{11})$ of $K3$ surfaces $X_t$ with a $\mu_{11}$ action;
see Oguiso-Zhang \cite{OZ-ord11}.
\end{rem}

\par
Now we give Kondo's examples of purely non-symplectic $\mu_I$ action on $K3$ surfaces with
$I$ a composite (see Kondo \cite{Ko92}, or Oguiso-Machida \cite{MO}, Prop. 4).

\begin{ex}\label{ex2}
Here are examples of pairs $(X_I, \mu_I)$ (the first few
bigger $I$ in Lemma \ref{index})
of the $K3$ surface $X_I$ with a $\mu_I = \langle g_I \rangle$ purely non-symplectic
action, given by its Weierstrass equation or as weighted hypersurface 
and the action on the coordinates; see Kondo \cite{Ko92}.

\par \vskip 1pc
$X_{66}: \,\, y^2 = x^3 + t(t^{11}-1)$, \,\, 
$g_{66}^*(x, y, t) = (\zeta_{66}^{40} x, \zeta_{66}^{27} y, \zeta_{66}^{54} t)$.

\par \vskip 1pc
$X_{50} = \{z^2 = x_0^6 + x_0x_1^5 + x_1 x_2^5\} \subset {\bold P}(1, 1, 1, 3)$,
\par 
\hskip 3pc $g_{50}^*[x_0:  x_1 :  x_2 : z] = [x_0 : \zeta_{25}^{20} x_1: \zeta_{25} x_2 : -z]$.

\par \vskip 1pc
$X_{44}: \,\, y^2 = x^3 + x +  t^{11}$, \,\, 
$g_{44}^*(x, y, t) = (\zeta_{44}^{22} x, \zeta_{44}^{11} y, \zeta_{44}^{34} t)$.

\par \vskip 1pc
$X_{54}: \,\, y^2 = x^3 + t(t^{9}-1)$, \,\, 
$g_{54}^*(x, y, t) = (\zeta_{27}^{2} x, -\zeta_{27}^{3} y, \zeta_{27}^{6} t)$.

\par \vskip 1pc
$X_{38}: \,\, y^2 = x^3 + t^7x + t$, \,\, 
$g_{38}^*(x, y, t) = (\zeta_{19}^{7} x, -\zeta_{19}y, \zeta_{19}^{2} t)$.

\par \vskip 1pc
$X_{48}: \,\, y^2 = x^3 + t(t^8-1)$, \,\, 
$g_{48}^*(x, y, t) = (\zeta_{48}^{2} x, \zeta_{48}^3y, \zeta_{48}^{6} t)$.
\end{ex}

\par
Indeed, these surfaces (and the action to some extent) are unique:

\begin{thm} (Machida-Oguiso \cite{MO}, Main Th $1$).
Suppose that $X$ is a $K3$ surface with a finite group $G \le \Aut(X)$
such that $I = I(X, G) \in \{66, 50, 44$, $33, 25\}$ (i.e., $\varphi(I) = 20$).

\par \noindent
$(1)$ The pair $(X, G)$ is isomorphic to the pair $(X_I, \langle g_I \rangle)$
in the example above if $I$ is even, and
$(X, G)$ is isomorphic to the pair $(X_{2I}, \langle g_{2I}^2 \rangle)$
if $I$ is odd.

\par \noindent
$(2)$ We have $\Aut(X) = G \cong \mu_I$ when $I$ is even, and
$\Aut(X) \cong \mu_{2I}$ when $I$ is odd.
\end{thm}

\begin{rem}
A slight weak result holds for $I = 54, 38, 48$ (with the unique surface $X_I$);
see Xiao \cite{Xi2}, Theorem, page 1.
\end{rem}

\par
So far, we have considered automorphisms of $K3$ surfaces with bigger order.
We now consider the opposite cases. We start with some examples.

\begin{ex}\label{ex3}
Let $\zeta = \zeta_4 = \sqrt{-1}$. Let $E = E_{\zeta} = {\bold C}/({\bold Z} + {\bold Z} \zeta)$
be the elliptic curve of period $\zeta$. Let $X_2 \rightarrow \overline{X}_2 
:= (E \times E)/ \diag(\zeta, -\zeta)$ be the minimal resolution of the quotient 
surface $\overline{X}_2$. Then $X_2$ is a $K3$ surface.
Let $g_2$ be the automorphism of $X_2$ induced by the action
$\diag(-1, 1)$ on $E \times E$. Then $g_2^* \omega_{X_2} = -\omega_{X_2}$
and the fixed locus (point wise) $X^{g_2}$ consists of 10 smooth rational
curves (see Oguiso-Zhang \cite{OZ-AJM}, Example 2 for details).
\end{ex}

\begin{ex}\label{ex4}
Let $\zeta = \zeta_3$. Let $E = E_{\zeta} = {\bold C}/({\bold Z} + {\bold Z} \zeta)$
be the elliptic curve of period $\zeta$. Let $X_3 \rightarrow \overline{X}_3 
:= (E \times E)/ \diag(\zeta, \zeta^2)$ be the minimal resolution of the quotient 
surface $\overline{X}_3$. Then $X_3$ is a $K3$ surface.
Let $g_3$ be the automorphism of $X_3$ induced by the action
$\diag(\zeta, 1)$ on $E \times E$. Then $g_3^* \omega_{X_3} = \zeta \omega_{X_3}$
and the fixed locus (point wise) $X^{g_3}$ consists of 6 smooth rational
curves and 9 isolated points (see Oguiso-Zhang \cite{OZ-AJM}, Example 1 for details).
\end{ex}

From the result above or intuitively, we see that a $K3$ surface $X$ with
a non-symplectic action of bigger order is uniquely determined by the  action.
The result below shows that the other extreme case may happen too.

\begin{thm} (Oguiso-Zhang \cite{OZ-AJM}, Theorem $4$).
Let $(X, g)$  be a pair of a $K3$ surface $X$ and an automorphism
$g$ of $X$ satisfying
(where the conditin $(3)$ is removable by Nikulin \cite{Ni-aut}, or Zhang \cite{Z-inv} Th $3$):

\par \noindent
$(1)$ $g^2 = id$,

\par \noindent
$(2)$ $g^* \omega_X = - \omega_X$,

\par \noindent
$(3)$ The fixed locus (point wise) $X^g$ consists of only (smooth) rational curves, and

\par \noindent
$(4)$ $X^g$ contains at least $10$ rational curves.

\par
Then $(X, g)$ is isomorphic to the pair $(X_2, g_2)$ up to isomorphisms.
\end{thm}

\begin{thm} (Oguiso-Zhang \cite{OZ-AJM}, Theorem $3$).
Let $(X, g)$  be a pair of a $K3$ surface $X$ and an automorphism
$g$ of $X$ satisfying:

\par \noindent
$(1)$ $g^3 = id$,

\par \noindent
$(2)$ $g^* \omega_X = \zeta_3 \omega_X$,

\par \noindent
$(3)$ The fixed locus (point wise) $X^g$ consists of only (smooth) rational curves
and possibly some isolated points, and

\par \noindent
$(4)$ $X^g$ contains at least $6$ rational curves.

\par
Then $(X, g)$ is isomorphic to the pair $(X_3, g_3)$ up to isomorphisms.
\end{thm}

\begin{rem}

\par \noindent
(1) The $K3$ surfaces $X_3$ and $X_2$ in the theorems above
have Picard number 20 and discriminants $|\det(\Pic(X_3))| = 3$ and
$|\det$ $(\Pic(X_2))| = 4$. Indeed, $X_3$ (resp. $X_2$) is the only
$K3$ surface with Picard number 20 and discriminant 3 (resp. 4).
They are called the {\it most algebraic $K3$ surfaces} in Vinberg \cite{Vi}.
For these two $X = X_i$, the infinite group $\Aut(X)$ has been calculated
by Vinberg.

\par \noindent
(2) A smooth $K3$ surface $X$ is called {\it singular} if the Picard number
$\rho(X) = 20$
(maximum). Such $X$ is uniquely determined by the transendental lattice
$T_X$; see Shioda-Inose \cite{SI}.
Moreover, it has been proved there that
a singular $K3$ surface has infinite group $\Aut(X)$.
See Oguiso's theorem below for a stronger result. 

\par \noindent
(3) Nikulin, Vinberg and Kondo  have determined all $K3$ surfaces $X$
with finite $\Aut(X)$; see Nikulin \cite{Ni-aut} and Kondo \cite{Ko89}.

\par \noindent
(4) In general, it is very difficult to calculate the full 
automorphism group $\Aut(X)$ for a $K3$ surface $X$.
This has been done by Keum \cite{Ke97} and Kondo \cite{Ko98} 
for a generic Jacobian Kummer surface
and by Keum-Kondo \cite{KK} for Kummer surfaces associated 
with products of some elliptic curves.

\par \noindent
(5) In their paper \cite{DK1}, Dolgachev and Keum have successfully
calculated the full automorphism group $\Aut(X)$ for general quartic Hessian surfaces,
by embedding the rank $16$ Picard lattice into the rank $26$ lattice 
$L$ of signature $(1,25)$ (the sum of the rank $24$ Leech lattice and the 
standard hyperbolic plane). See also Kondo \cite{Ko-max} for the employment of
such $L$ for Kummer surfaces.
\end{rem}

We end this section with the following striking result
in (Oguiso \cite{Og-free}); see the same paper for more 
results in the new directions.

\begin{thm} (Oguiso \cite{Og-free}).
Let $X$ be a singular $K3$ surface and let $H := Hilb^n(X)$.
Then each of $\Aut(X)$, $\Aut(H)$ and $\Bir(M)$ contains
the non-commutative free group ${\bold Z} * {\bold Z}$.
\end{thm}

\section{\bf $K3$ groups of mixed type}

So far, we have considered $K3$ groups of purely symplectic or purely non-symplectic type.
In this section, we shall consider $K3$ groups $G$ of the mixed type, i.e., the case where
neither $G_N$ nor $\mu_I$ (with $I = I(X, G)$) in the basic exact sequence of the Introduction
becomes trivial. On the one hand, $|G_N| \le |M_{20}| = 960$ by Theorem \ref{Mu-sym}.
On the other hand, $|\mu_I| \le 66$ by Lemma \ref{index}.
A finite $K3$ group $G$ (and the pair $(X, G)$) 
of the largest order has been determined by Kondo:

\begin{thm} (Kondo \cite{Ko-max}).
Let $X$ be a $K3$ surface and $G \le \Aut(X)$ a finite subgroup.

\par \noindent
$(1)$ We have $|G| \le 4 \cdot 960$.

\par \noindent
$(2)$ Suppose that $|G| = 4 \cdot 960$. Then $X$ equals the Kummer surface
$Km(E_{\sqrt{-1}} \times E_{\sqrt{-1}})$ where 
$E_{\sqrt -1} = {\bold C}/({\bold Z} + {\bold Z}\sqrt{-1})$, and
$G = M_{20} . \mu_4$. Moreover, the group $G$ and the action of $G$ on $X$
are unique up to isomorphisms.
\end{thm}

Besides the fact that $M_{20}$ is the finite symplectic $K3$ group of the largest order,
it is also the the finite perfect $K3$ group of the largest order.
Here are all of the finite perfect $K3$ groups (cf. Xiao \cite{Xi1}, the list):
$$M_{20}, \hskip 0.5pc A_6, \hskip 0.5pc L_2(7) = \text{\rm PSL}_2(7), \hskip 0.5pc A_5.$$
The first three groups are among the 11 maximum ones in Theorem \cite{Mu-sym}.
The last three groups above are also the only non-abelian finite simple $K3$ groups.
The last two groups happen to be the two smallest (in terms of order) 
non-abelian finite simple groups.
These four groups (together with a bigger $\mu_I$)
determine the surface $X$ uniquely (see the theorems below).

Let us start with an example of $K3$ group of mixed type.

\begin{ex}\label{ex5}
Let $C_{168} = \{x_1x_2^3 + x_2x_3^3 + x_3x_1^3 = 0\} \subset {\bold P}^2$
be the Klein quartic curve of genus $g = 3$. A well known theorem
says that $|\Aut(C)| \le 84(g-1)$ for any curve $C$ of genus $g \ge 2$.
This Klein quartic attains the maximum $|\Aut(C_{168})| = 84 (3 - 1) = 168
= |L_2(7)|$ (see Oguiso-Zhang \cite{OZ-168} for the action of $L_2(7)$ on $C_{168}$).
Consider the quartic $K3$ surface 
$X_{168} = \{x_0^4 + x_1x_2^3 + x_2x_3^3 + x_3x_1^3 = 0\} \subset {\bold P}^3$.
This is also the Galois ${\bold Z}/(4)$-cover of ${\bold P}^2$ branched
along the Klein quartic $C_{168}$. Clearly,
there is a faithful $L_2(7) \times \mu_4$ action on $X_{168}$.
We shall call the pair $(X_{168}, L_2(7) \times \mu_4)$ the Klein-Mukai pair.
\end{ex}

\begin{thm} (Oguiso-Zhang \cite{OZ-168}, Main Th, Prop. $1$).
Let $X$ be a $K3$ surface such that $G : = L_2(7)\le \Aut(X)$.
Let $\widetilde{G} \le \Aut(X)$ be a finite subgroup containing $G$.

\par \noindent
$(1)$ $G$ is normal in $\widetilde{G}$ with $|{\widetilde G}/G| \le 4$.

\par \noindent
$(2)$ Suppose that $|\widetilde{G}/G| = 4$. Then $(X, \widetilde{G})$
is isomorphic to the Klein-Mukai pair $(X_{168}, L_2(7) \times \mu_4)$.
\end{thm}

\par
Next we turn to the another simple group $A_6$, also called the Valentiner's group
according to Dolgachev and Gizatullin
(see the reference in Keum-Oguiso-Zhang \cite{KOZ2}).
It is a very important group and also the root of troubles in the theory of simple groups.
Indeed, the $A_6$ is at the junction of the three simple groups series:
recall that $A_6 = [M_{10}, M_{10}] = \text{\rm PSL}(2, 9)$.

\begin{thm} (Keum-Oguiso-Zhang \cite{KOZ1}, Th $3.1$, Prop $4.1$, Th $5.1$).
$(1)$ Let $F$ be the unique $K3$ surface of Picard number $20$ whose
transcendental lattice $T_F$ has the intersection matrix $\diag[6, 6]$
(see Shioda-Inose \cite{SI}).
Then there is a group $\widetilde{A}_6 := A_6 : \mu_4$
(a split extension of $A_6$ by $\mu_4$) and a faithful action of it
on $F$ (see Keum-Oguiso-Zhang \cite{KOZ1} Def $2.7$, Th $3.1$ or \cite{KOZ2} Th $3.2$
for the action and the group structure).

\par \noindent
$(2)$ Suppose that $\widetilde{G}$ is a finite $K3$ group 
containing $A_6$. Then $A_6$ is normal in $\widetilde{G}$,
and the group $\widetilde{G}/A_6$ is cyclic of order $\le 4$.

\par \noindent
$(3)$ Suppose that $(X, \widetilde{G})$ is a pair of a $K3$ surface
and a subgroup $\widetilde{G} \le \Aut(X)$ containing $A_6$ as a subgroup
of maximum index $4$. Then $(X, \widetilde{G})$ is isomorphic
to the pair $(F, \widetilde{A}_6)$ in $(1)$.
\end{thm}

\begin{rem}
The $F$ in the theorem above is also the minimal resolution of the hypersurface
in ${\bold P}^1 \times {\bold P}^2$ (with coordinates $([S:T], [X:Y:Z])$)
given by the equation (see Keum-Oguiso-Zhang \cite{KOZ1} Prop 3.5):
$$S^2(X^3+Y^3+Z^3) - 3(S^2+T^2)XYZ = 0.$$

The action of $\widetilde{A}_6$ is still invisible from the equation
(due to the trouble nature of $A_6$?)
It is indeed constructed lattice-theoretically, where the Leech lattice
and some deep results of S. Kondo \cite{Ko-lat} and R. E. Borcherds \cite{Bo} are utilized.
\end{rem}

\par
In the case of the smallest non-abelian simple group $A_5$, 
we have the complete classification of its finite $K3$ over-groups:

\begin{thm} (Zhang \cite{Z-A5} Th A, or \cite{Z-dol} Th E).
Suppose that $G$ is a finite $K3$ group containing $A_5$ as a normal subgroup.
Then $G$ equals one of the following
(for their geometric realization, see Zhang \cite{Z-A5} Example $1.10$):
$$A_5, \hskip 0.5pc S_5, \hskip 0.5pc A_5 \times \mu_2, \hskip 0.5pc S_5 \times \mu_2.$$
\end{thm}

We now give an example of the largest (in terms of the order) finite solvable $K3$ group.

\begin{ex}\label{ex6}
Let $X_4 := \{x_1^4 + x_2^4 + x_3^4 + x_4^4 = 0 \} \subset {\bold P}^3$
be the Fermat quartic $K3$ surface.
Then $(\mu_4)^4$ acts diagonally on the coordinates while $S_4$
permutes the coordinates. So there is a natural
action of $\widetilde{F}_{384} = (\mu_4^4 : S_4) / \mu_4 = (\mu_4^4 / \mu_4) : S_4$
on $X_4$ (see Mukai \cite{Mu-sym}).
This $\widetilde{F}_{384}$ is a solvable group of order $4^3 \cdot 4!  = 2^9 \cdot 3 = 384 \cdot 4$. 
Note that the $2$-Sylow subgroup $\widetilde{F}_{128}$ of $\widetilde{F}_{384}$
is of course nilpotent and of order $2^9 = 128 \cdot 4$.
Let $G = \widetilde{F}_{128}$, or $\widetilde{F}_{384}$, then
in notation of the basic exact sequence in the Introduction,
we have $\mu_{I(X, G)} = \mu_4$ and $G_N = F_{128}$, or $F_{384}$ respectively 
(as defined in Mukai \cite{Mu-sym};
note also that $F_{384}$ is also the second largest among the 11 maximum ones in 
Theorem \ref{Mu-sym}). For both $G$, one can show that 
$\Pic(X_4)^G = {\bold Z} H$ and $H^2 = 4$ (see Oguiso \cite{Og-comp}, Th 1.2).
\end{ex}

Now we state a pretty result due to Oguiso:

\begin{thm} (Oguiso \cite{Og-comp}, Th $1.2$).
\par \noindent
$(1)$ Let $X$ be a $K3$ surface and $G \le \Aut(X)$ a finite solvable subgroup.
Then $|G| \le 2^9 \cdot 3$. Moreover, if $|G| = 2^9 \cdot 3$
then $(X, G)$ is isomorphic to the pair $(X_4, \widetilde{F}_{384})$
(with the standard action) in the example above.

\par \noindent
$(2)$ Let $X$ be a $K3$ surface and $G \le \Aut(X)$ a finite nilpotent subgroup.
Then $|G| \le 2^9$. Moreover, if $|G| = 2^9$
then $(X, G)$ is isomorphic to the pair $(X_4, \widetilde{F}_{128})$
(with the standard action) in the example above.
\end{thm}

\section{\bf Log Enriques surfaces}

The concept of Log Enriques surface was first introduced in
Zhang \cite{Z-enriq}, Def 1.1. In the smooth case, they are
just abelian surface (index 1), $K3$ surface (index 2), Enriques surface (index 2)
and hyperelliptic surfaces (index $2, 3, 4$, or $6$).
Log Enriques surface is closely related to the
study of purely non-symplectic $K3$ groups.

\begin{definition}
A normal projective surface $Y$ with at worst quotient singularities is called
a {\it log Enriques surface} if the irregularity $h^1(Y, \OO_Y)$ $= 0$
and if a positive multiple of the canonical divisor $m K_Y \sim 0$ (linearly equivalence).
$I = I(Y) = \min \{m \in {\bold Z}_{>0} \, | \, m K_Y \sim 0\}$
is called the {\it index} of $Y$.

There is a Galois ${\bold Z}/(I)$-cover
$\pi : X = Spec \oplus_{i=0}^{I-1} \OO_Y(-iK_Y) \rightarrow Y$
which is unramified outside non-Du Val singular locus of $Y$,
i.e., outside the set of singular points which is not rational double.
This $X$ is called the {\it canonical cover} of $Y$. One sees
that $g^* \omega_X = \eta_I \omega_X$, where $\Gal(X/Y) = \langle g \rangle \cong \mu_I$
and $\eta_I$ is a primitive $I$-th root of $1$.
\end{definition}

\begin{rem}
\par \noindent
(1) In general, the canonical cover $X$ is a normal surface with at worst
Du Val singularities and with $K_X \sim 0$.
So either $X$ is a (smooth) abelian surface, or a $K3$ surface
with at worst Du Val singularities.

\par \noindent
(2) Log Enriques surfaces are degenerate fibres
of families of Kodaira dimension 0. Also the base space of elliptically
fibred Calabi-Yau threefold $\Phi_D : Z \rightarrow Y$ with $D . c_2(Z) = 0$
are necessarily log Enriques surfaces (see Oguiso \cite{Og-IJM}).

\par \noindent
(3) When the canonical cover $X$ of $Y$ is abelian, one has $I(Y) = 3$, or $5$
and such $Y$ is unique up to isomorphisms; for examples or proof, see
Zhang \cite{Z-enriq} Example 4.2 and Theorem 4.1, or 
Blache \cite{Bl}, (1.2) and Theorem C.

\par \noindent
(4) In general, one has $I(Y) \le 21$; see Blache \cite{Bl} Th C, and also 
Zhang \cite{Z-enriq} Lemma 2.3. If $I(Y)$ is prime, then $I(Y) \le 19$
(Lemma \ref{index}).
Examples of prime $I(Y)$ are given in Zhang \cite{Z-enriq} and Blache \cite{Bl}.
When $I(Y)$ is prime, the classification of $\Sing(X)$
for the canonical cover $X$ of $Y$, is done in Zhang \cite{Z-enriq2}, Main Th.

\par \noindent
(5) When the canonical cover $X$ of $Y$ is a normal $K3$, we 
let $\widetilde{X} \rightarrow X$ be the minimal resolution.
Then the determination of $Y$ is almost equivalent to that of $(\widetilde{X}, \mu_I)$,
where $\mu_I = \Gal(X/Y)$ acts purely non-symplectically on the (smooth) $K3$
surface $\widetilde{X}$. Clearly, $20 \le \rho(\widetilde{X}) = r(Y) + \rho(X) \ge r(Y) + 1$
with $r(Y)$ the number of components in the exceptional divisor of the resolution
$\widetilde{X} \rightarrow X$.
Here $\rho(Z)$ denotes the Picard number. Hence $r \le 19$.

\par \noindent
(6) For rational log Enriques surfaces $Y$ with maximum $r(Y) = 19$,
the type $\Sing(X)$ has been determined. They are one of 7 Dynkin types:
$$D_{19}, \,\, D_{16} + A_3, \,\, D_{13} + A_6, \,\, D_7 + A_{12}, \,\, 
D_7 + D_{12}, \,\, D_4 + A_{15}, \,\, A_{19}.$$
For each of the 7 Dynkin types above, there is exactly one (up to isomorphisms)
rational log Enriques surface $Y$ such that the singular locus $\Sing(X)$
of the canonical cover $X$ of $Y$ has that Dynkin type;
see Oguiso-Zhang \cite{OZ-AJM} Th 1, Th 2 and \cite{OZ-mathZ} Main Th.
\end{rem}

Start with a log Enriques surface $Y$, by passing to a maximum crepant
partial resolution (i.e., resolving Du Val singularities), we may assume
the maximality condition for $Y$:

\par \vskip 1pc \noindent
(*) \,\, Any birational morphism $Y' \rightarrow Y$ from another
log Enriques surface $Y'$ must be an isomorphism.

\par \vskip 1pc
We now state a uniqueness result for the three largest prime indices:

\begin{thm} (Oguiso-Zhang \cite{OZ-PAMS}, Cor. 4).
Let $Y$ be a log Enriques surface with $I = I(Y) = 19, 17$ or $13$ 
and satisfying the maximality condition (*) above.

\par \vskip 1pc \noindent
$(1)$ We have $Y \cong X_{I}/\langle g_I \rangle$ when $I = 19$ or $17$.

\par \vskip 1pc \noindent
$(2)$ We have $Y \cong \overline{X}_{13}/\langle g_{13} \rangle$ when $I = 13$.

\par \vskip 1pc \noindent
Here the pairs $(X, \langle g_I \rangle)$ are given in Example \ref{ex1},
and $\overline{X}_{13}$ is obtained from
${X}_{13}$ by contracting the unique rational curve in the
fixed locus (point wise) $X_{13}^{g_{13}}$.
\newline
(see Zhang \cite{Z-enriq} Ex $5.7 - 5.8$, \cite{Z-enriq2} Ex $7.3$
for different constructions).
\end{thm}

We end the paper with the final remark.

\begin{rem}

\par \noindent
(1) In their paper \cite{DK}, Dolgachev-Keum have extended Mukai's classification
Theorem \ref{Mu-sym}
to characteristic $p$ case. For smaller $p$, some new groups of automorphisms (not in
Mukai's list) appear, see ibid and Dolgachev-Kondo \cite{DKo}.

\par \noindent
(2) In his paper \cite{Og-JAG}, Theorem 1.5, Oguiso considers the subtle behavior of
$\Aut(X_t)$ for $K3$ surfaces in a 1-dimensional family.
On the one hand, for generic $t$ the group $\Aut(X_t)$ contains a fixed subgroup
with bounded index. On the other hand, he also produces an example where $\Aut(X_t)$
is infinite for generic $t$ but is finite for some special $t$. 

\par \noindent
(3) Many interesting results could not be included in the paper due to the
limitation of the knowledge of the author and the
constraint of the space.
\end{rem}
%%%%%%%%%%%%%%%%%%%%%%%%%%%%%%%%%%%%%%%%%%%%%%%%%%%%%%%

\end{document}